\documentclass[10pt,a4paper]{article}

\usepackage{epsfig}
\usepackage{amsmath}
\usepackage{amsfonts}
\usepackage{latexsym,amsthm}

\DeclareSymbolFont{EulerScript}{U}{eus}{m}{n}
\DeclareSymbolFontAlphabet\mathscr{EulerScript}

\begin{document}

\newcommand{\rk}{\mathsf{rk}}
\newcommand{\RR}{\mathsf{R}}
\newcommand{\LL}{\mathsf{L}}
\newcommand{\N}{\mathbb N}
\newcommand{\R}{\mathbb R}
\newcommand{\Z}{\mathbb Z}
\newcommand{\Q}{\mathbb Q}
\newcommand{\C}{\mathbb C}

\newcommand{\esp}{\vskip .3cm \noindent}

\newcommand{\el}{{\cal L}}
\newcommand{\de}{{\cal D}}
\newcommand\Cfty{${\cal C}^{\infty}$}
\newcommand{\emm}{{\cal M}}
\newcommand{\be}{{\cal B}}
\newcommand{\ess}{{\cal S}}
\newcommand{\enn}{{\cal N}}
\newcommand{\ka}{{\cal K}}
\newcommand{\emmh}{\hat {\cal M}}
\newcommand{\RD}{\text{\rm RD}}
\newcommand{\definition}{\overset{\text{\rm def}}{=}}
\newcommand{\graph}{\text{\rm graph}}

\def\CC#1{${\cal C}^{#1}$}
\def\wt#1{\widetilde{#1}}
\def\wh#1{\widehat{#1}}
\def\verti#1#2{\scriptsize\binom{#2}{#1}}

\newtheorem{prop}{Proposition}[section]
\newtheorem{lemma}[prop]{Lemma}
\newtheorem{cor}[prop]{Corollary}
\newtheorem{defi}[prop]{Definition}
\newtheorem{thm}[prop]{Theorem}

\newtheoremstyle{mypl}{10pt}{10pt}{\it}{}{\em}{.}{ }{}
\theoremstyle{mypl}
\newtheorem{claim}{Claim}[prop]
\renewcommand\theclaim{\Alph{claim}}
\newtheorem*{claim*}{Claim}

\newtheoremstyle{myrem}{10pt}{10pt}{\rm}{}{\bf}{.}{ }{}
\theoremstyle{myrem}
\newtheorem*{rem}{Remark}

\renewcommand{\theprop}{\thesection.\arabic{prop}}

\centerline{\bf LONG LINE KNOTS}
\baselineskip=13pt
\vspace*{0.37truein}
\centerline{\footnotesize MATHIEU BAILLIF \& DAVID CIMASONI}
\vspace*{0.37truein}
\begin{abstract}
  We study continuous embeddings of the long line $\LL$ into $\LL^n$ $(n\ge 2)$ up to ambient isotopy of $\LL^n$. We define the direction of an embedding
  and show that it is (almost) 
  a complete invariant in the case $n=2$ for continuous embeddings, and in the case $n\ge 4$ for differentiable ones. Finally, we
  prove that the classification of smooth embeddings $\LL\to\LL^3$ is equivalent to the classification of classical oriented knots. 
\end{abstract}
\section*{Introduction}

%%%% intro

Consider the following general problem:
\begin{equation}\tag{$\star$}
  \begin{array}{c}
    \text{Given finite dimensional manifolds $X,Y$, classify}\\
    \text{the embeddings $X\to Y$ up to ambient isotopy of $Y$.}
  \end{array}
  \label{problem}
\end{equation}
An instance of this problem is classical
knot theory, where $X=\mathbb{S}^1$, $Y=\mathbb{S}^3$ and the embeddings are assumed smooth (or PL).
Another instance is higher dimensional knot theory, where
$X=\mathbb{S}^k$ and $Y=\mathbb{S}^n$.
These fields have been very popular among mathematicians for over a century; there is a considerable litterature
about the classification of (topological, differentiable or PL) embeddings from a $k$-sphere to an $n$-sphere.

On the other hand, the study of this problem for non-metrizable manifolds remains to be done. 
Even the study of homotopy and isotopy classes of maps of non-metrizable
manifolds seems to be at its very beginning. To our knowledge, David Gauld \cite{Gauld} was the first ever to publish a paper on the subject
(see also \cite{DeoGauld}). He investigated homotopy classes of maps of the {\em long line} $\LL$ (and of
the {\em long ray} $\RR$) into itself (see Definition \ref{Longlinedef} below), showing that
there are exactly $2$ such classes for the long ray, and $9$ for the long line. 
He also proved that all embeddings $\RR\to\RR$ are isotopic, and that there are 
$2$ isotopy classes
of embeddings $\LL\to\LL$. These results provide a solution to Problem (\ref{problem}) for $X=Y=\RR,\LL$.

In this paper, we investigate Problem (\ref{problem}) with $X=\LL$, $Y=\LL^n$ and $X=\RR$, $Y=\RR^n$.
We introduce a numerical invariant of embeddings $\RR\to\RR^n$ or $\LL\to\LL^n$ called the {\it direction},
see Definition \ref{def direction} and Theorem \ref{thm:invariant}.
Roughly speaking, two knots (that is, two embeddings) have the same direction if for all $1\le i\le n$, their projections on 
the $i$-th coordinate 
are either both cofinal or both bounded.

Our first results are that the direction is (almost) a complete invariant in the cases
$n=2$ for continuous embeddings (Theorems \ref{thm1} and \ref{corthm1}) and
$n\ge 4$ for smooth embeddings (Theorem \ref{thm2} and \ref{corthm2}). It follows that for these $n$, 
there are exactly $2^n-1$ long ray knots and $(3^n-1)^2+(n-2)2^{n-1}$ long line knots (in the classes specified above).
We also prove that there are $7$ smooth knots $\RR\to\RR^3$ (Theorem \ref{thm3}).
Finally, our last result shows that the classification of differentiable embeddings $\LL\to\LL^3$ reduces to classical oriented knot theory
(see Proposition \ref{prop:ass} and Theorem \ref{thm4}).

The paper is organized as follows: Section \ref{secdef} deals with the definition of some basic objects, such as the long ray $\RR$,
the long line $\LL$, and the equivalence relation for knots. In Section \ref{sectools}, we prove several technical lemmas 
(mainly partition and covering properties) that are used in
Section \ref{secres}, which contains all the main results.

The authors wish to acknowledge Claude Weber, David Gauld and Ren\'e Binam\'e.

% section 1

\section{Definitions}\label{secdef}
For any ordinal $\alpha$, let us denote by $W(\alpha)$ the set of ordinals strictly smaller than $\alpha$.
\begin{defi}
\label{Longlinedef}
  The (closed) long ray $\RR$ is the set $W(\omega_1)\times[0,1[$,
  equipped with the lexicographic order and the order topology.
  The long line $\LL$ is the union of two copies $\LL_-,\LL_+$ of $\RR$ glued at $(0,0)$.
  We put the reverse order on $\LL_-$, so that $\LL$ is totally ordered.
\end{defi}
We will often identify $\alpha\in W(\omega_1)$ with $(\alpha,0)\in\RR$ and similarly for $\LL$.
Recall that $\RR$ and $\LL$ are non-metrizable, non-contractible and sequentially compact.
Also, $\LL$ and $\RR$ can be given a structure of oriented \CC{\infty} manifold. We will assume
throughout the text that we are given a {\em fixed} maximal atlas $\{U_j,\psi_j\}$ on $\LL$ (and on $\RR$)
with $U_j\ni 0$ for all $j$.
The atlas on $\LL^n$ (or $\RR^n$) is then assumed to be $\{(U_j)^n,(\psi_j,\dots,\psi_j)\}$,
the so-called `$n$-th power structure'.\footnote{These precisions are necessary, since P.J. Nyikos \cite{Nyikos:1992} 
showed  that there are
uncountably many non-equivalent differentiable structures on $\LL$ (and thus, on $\LL^n$). Moreover, it is not clear that any differentiable structure on
$\LL^n$ is equivalent to a product of structures on $\LL$.}
This assumption is, in our opinion, quite natural.
For instance, it
ensures that the maps $\RR\to\RR^2$, given by
$x\mapsto (x,x)$ and $x\mapsto (x,c)$ are smooth, where $c$ is any constant.
Throughout the text, $\pi_i$ will denote the projection on the $i$-th coordinate of $\LL^n$ or $\RR^n$.
\begin{defi}
  We call an embedding from $X$ to $Y$ a $(Y,X)$-knot.
  Two $(\LL^n,\LL)$-knots $f,g$ are equivalent, which we denote by $f\sim g$,
  if there is an isotopy $\phi_t$ $(t\in[0,1])$ of $\LL^n$
  such that $\phi_0=id$ and $\phi_1\circ f=g$.
\end{defi}
Since $\RR$ is a manifold with boundary, we have to be more careful with the definition of equivalent
$(\RR^n,\RR)$-knots.
One possibility is to consider only embeddings of $\RR$ in the interior of $\RR^n$. Here is another
way to state the same equivalence relation:
\begin{defi}
  Let $\RR'$ be the set $(\{-1\}\times [0,1[)\sqcup\RR$ equipped with the order topology.
  Two $(\RR^n,\RR)$-knots $f,g$ are equivalent, which we denote by $f\sim g$,
  if there is an isotopy $\phi_t$ $(t\in[0,1])$ of $(\RR')^n$
  such that $\phi_0=id$ and $\phi_1\circ f=g$.
\end{defi}

%%%%%% sec 2

\section{Tools}\label{sectools}

Let us begin by recalling the following well-known lemma; the proof of the corresponding statement
for ordinals can be found in any book on set theory.
\begin{lemma}\label{closedcofinal}
  Let $\{E_m\}_{m<\omega}$ be closed and cofinal subsets of $\RR$. Then $\bigcap_{m<\omega}E_m$ is also closed and cofinal in $\RR$.
  \hfill $\Box$
\end{lemma}
\noindent
We shall now investigate partition and covering properties of embeddings $\RR\to\RR^n$.
\begin{lemma}\label{partition}
  Let $\{f_k\}_{k\in K}$ be a finite or countable family of continuous maps $\RR\to\RR$. \\
  a) If each $f_k$ is bounded, there
  is a $z$ in $\RR$ such that $f_k$ is constant on $[z,\omega_1[$ for all $k\in K$. \\
  b) If each $f_k$ is cofinal, there is a
  cover $\mathcal{P}=\{[x_\alpha,x_{\alpha+1}]\}_{\alpha<\omega_1}$ of $\RR $ such that
  $x_\alpha<x_\beta$ if $\alpha<\beta$, 
  $x_\beta =\sup_{\alpha<\beta} x_\alpha$ if $\beta$ is a limit ordinal, 
  and for all $k\in K$,
  $f_k([x_\alpha,x_{\alpha+1}])=[x_\alpha,x_{\alpha+1}]$ for $\alpha>0$
  and $f_k([0,x_1])\subset [0,x_1]$.
\end{lemma}
\noindent{\em Proof:}
a) By \cite[Lemma 3.4 (iii)]{Nyikos:1984}, there exists $z_k\in \RR$ such that $f_k$ is constant on $[z_k,\omega_1[$. Take $z=\sup_{k\in K}z_k$.\\
b) This is a consequence of the following claims:
\begin{claim}\label{clm1}
  $A_k\definition\{ x\in \RR \, |\,  f_k(y)\le x \,\,\forall y\le x\}$ is closed and
  cofinal in $\RR $.
\end{claim}
\begin{claim}\label{clm2}
 $B_k\definition\{ x\in \RR \, |\,  f_k(y)\ge x \,\,\forall y\ge x\}$ is closed and cofinal
 in $\RR $.
\end{claim}
Indeed, these two claims together with Lemma \ref{closedcofinal} imply that the set
$E=\bigcap_{k\in K}(A_k\cap B_k)=\{x\in\RR\,|\,f_k([0,x])\subset[0,x]\mbox{ and }f_k([x,\omega_1[)\subset[x,\omega_1[\,\, \forall k\in K\}$ 
is closed and cofinal. We then define by
transfinite induction
$x_\alpha\in E$ for each $\alpha<\omega_1$ as follows.
Set $x_0=0$, 
$x_{\alpha+1}=\min (E\cap [x_\alpha+1,\omega_1[)$,
and
$x_\beta=\sup_{\alpha<\beta}x_\alpha$ if $\beta$ is a limit ordinal.
(Recall that for any limit ordinal $\beta<\omega_1$, there is a sequence $\{\alpha_i\}_{i<\omega}$, with $\alpha_i<\beta$, such
that $\lim_{i\to\infty}\alpha_i=\beta$. Therefore, $x_\beta=\lim_{i\to\infty}x_{\alpha_i}$ belongs to $E$.) 
\esp
{\it Proof of Claim \ref{clm1}:} Closeness is obvious (recall that a sequentially closed subset of $\RR$ is closed).
Given $z$ in $\RR$, one can choose a sequence $\{x_j\}_{j<\omega}$
such that $z < x_j \le x_{j+1}$  and
$\sup_{[0,x_j]} f_k \le x_{j+1}$. The sequence converges to a point $x>z$,
and $x\in A_k$.
\esp
{\it Proof of Claim \ref{clm2}:} Again, closeness is obvious. Fix $z$ in $\RR$;
for any $x\in\RR$, there is a $y\in\RR $ with $y\ge x$ such that
$f_k([y,\omega_1[)\subset [x,\omega_1[$.
(Otherwise, the set of $y$ such that $f_k(y)< x$ is cofinal and closed;
by
cofinality of $f_k$, so is
the set of $y'$ such that $f_k(y')> x$. By Lemma \ref{closedcofinal},
there would exist $v$ such that $f_k(v)< x$ and $f_k(v)> x$, which is impossible.)
We then define sequences $x_j,y_j$ ($j<\omega$) such that
$z < x_j\le y_j\le x_{j+1}$ and $f_k([y_j,\omega_1[)\subset [x_j,\omega_1[$.
The sequences $x_j,y_j$ converge to the same point $x>z$, which belongs to $B_k$.
\hfill $\Box$
\begin{lemma}\label{pearl necklace lemma}
  Let $\mathcal{P}=\{[x_\alpha,x_{\alpha+1}]\}_{\alpha<\omega_1}$ be a cover of $\RR$ with
  $x_\alpha<x_\beta$ if $\alpha<\beta$ and $x_\beta=\sup_{\alpha<\beta} x_\alpha$ if $\beta$ is a limit ordinal,
  and let $S_\alpha$ be the cube $[x_\alpha,x_{\alpha+1}]^n$.
  Then, for all $0<\alpha<\omega_1$, there are
  closed subsets $B_\alpha$ of $\RR^n$, diffeomorphic to a compact ball in $\R^n$,
  such that 
  \begin{itemize}
    \item[$\bullet$]  
    $B_\alpha\supset S_\alpha$,
    \item[$\bullet$]
    $\partial B_\alpha\cap S_\alpha=\{(x_\alpha,\dots,x_\alpha),(x_{\alpha+1},\dots,x_{\alpha+1})\}$,
    \item[$\bullet$]$\partial B_\alpha\cap S_\beta=S_\alpha\cap S_\beta$ if $\alpha\not=\beta$.
  \end{itemize}
\end{lemma}
\begin{figure}[h]
\begin{center}
  \epsfig{figure=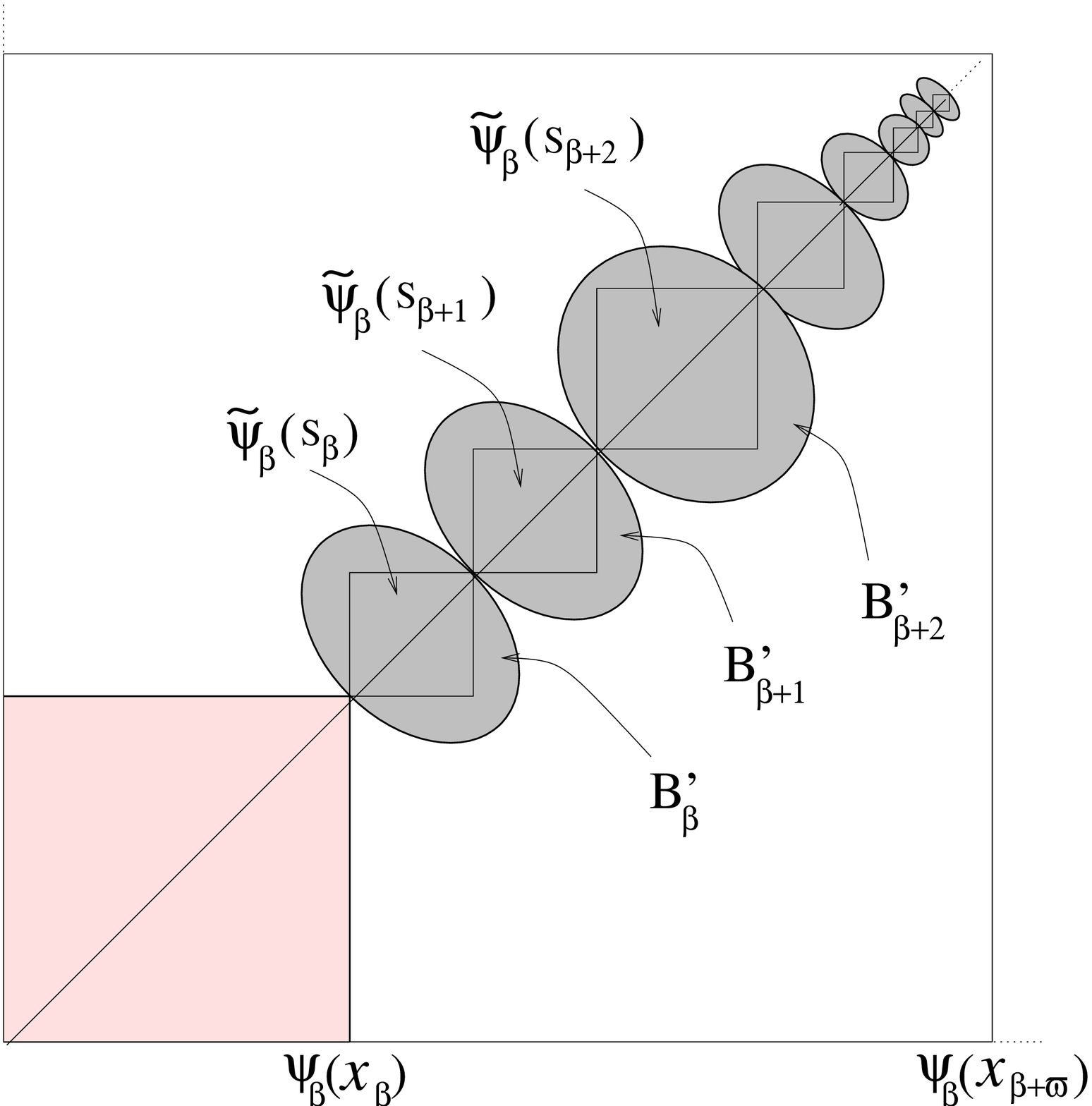,  height=5cm}
  \caption{Lemma \ref{pearl necklace lemma} for $n=2$.}
 \label{fig:pearl}
\end{center}
\end{figure}
\noindent{\em Proof:}
Let $\beta$ be a limit ordinal or $1$. Consider $\psi_\beta:U_\beta\to\R$ the smallest chart such that
$U_\beta\supset[0,x_{\beta+\omega}]$ (such a chart always exists), and
$\wt{\psi}_\beta=(\psi_\beta,\dots,\psi_\beta):(U_\beta)^n\to\R^n$, which belongs
to the atlas of $\RR^n$.
Clearly, $\wt{\psi}_\beta(S_\alpha)$ is a cube in $\R^n$ for any $\alpha<\beta+\omega$.
We can therefore form ellipsoids $B_\alpha '$ for $\beta\le\alpha<\beta+\omega$ as in Figure \ref{fig:pearl},
and define $B_\alpha=\wt{\psi}_\beta^{-1}(B_\alpha ')$ for $\beta\le\alpha<\beta+\omega$.
Note that $B_\beta '\cap \wt{\psi}_\beta([0,x_\beta]^n)=\{\wt{\psi}_\beta(x_\beta)\}$.
Since $W(\omega_1)\backslash\{0\}=[1,\omega[\,\sqcup\,\,\bigsqcup_{\beta \text{ lim.}}[\beta,\beta+\omega[$,
we obtain the desired properties.
\hfill $\Box$
\begin{figure}[h]
\begin{center}
  \epsfig{figure=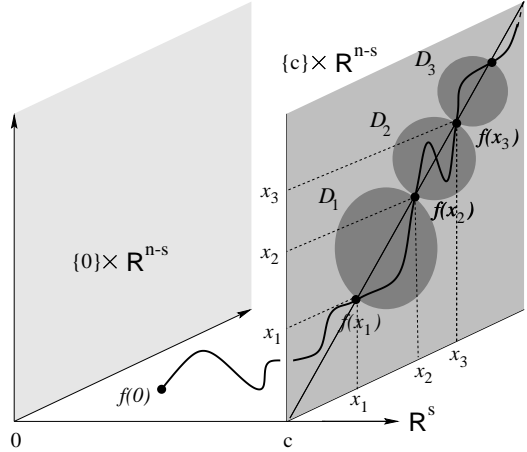,  height=6cm}
  \caption{Lemma \ref{boxing lemma} for $n=3$, $s=1$.}
  \label{fig:boxing}
\end{center}
\end{figure}
\begin{lemma}\label{boxing lemma}
Let $f:\RR\to\RR^n$ be continuous, with $\pi_i\circ f$ bounded for $i=1,\dots,s$ $(s<n)$ and
cofinal for $i=s+1,\dots,n$. Then, there exists a cover
$\mathcal{P}=\{[x_\alpha,x_{\alpha+1}]\}_{\alpha<\omega_1}$ of $\RR$ and
a cover
$\mathcal{D}=\{D_\alpha\}_{0<\alpha<\omega_1}$ of $f([x_1,\omega_1[)$ with
$D_\alpha$ diffeomorphic to the compact ball in $\R^n$,  
such that (see Figure \ref{fig:boxing}):
\begin{itemize}
\item[$\bullet$]
  $f(x)\subset\{c\}\times\RR^{n-s}\;\forall x\ge x_1$, for some fixed $c\in\RR^s$,
\item[$\bullet$]  
  $f(x_\alpha)=(c,x_\alpha,\dots,x_\alpha)$ $\forall \alpha>0$, 
%\item[$\bullet$]
%  $\pi_i\circ f([x_\alpha,x_{\alpha+1}])=[x_\alpha,x_{\alpha+1}]$ $\forall \alpha>0$
%  and $\pi_i\circ f([0,x_1])\subset [0,x_1]$, for $i=s+1,\dots,n$,
%\item[$\bullet$]
%  $f([0,x_1[)\cap D_\alpha$  is empty $\forall\alpha>0$,
\item[$\bullet$]
  $f^{-1}(D_\alpha)=[x_\alpha,x_{\alpha+1}]$ and $f^{-1}(\partial D_\alpha) = \{x_\alpha,x_{\alpha+1}\}$ $\forall\alpha>0$,
 
\item[$\bullet$]
  $D_\alpha\cap D_\beta=\left\{\begin{array}{ll}
            \emptyset &\text{\rm if }\alpha\not=\beta\pm 1\\
            (c,x_\alpha,\dots,x_\alpha)&\text{\rm if }\alpha=\beta+ 1.
            \end{array}
            \right.$
\end{itemize}
\end{lemma}
\noindent{\em Proof:}
By Lemma \ref{partition} a), there is some $z\in\RR$ such that $\pi_i\circ f$ is
constant on $[z,\omega_1[$ for $i=1,\dots,s$. Set $c_i=\pi_i\circ f(z)$ and $c=(c_1,\dots,c_s)$.
For $i=s+1,\dots,n$,
take the cover $\mathcal{P}$ for the family $\{\pi_i\circ f\}_{i=s+1,\dots,n}$ with $x_1>z$ given by
Lemma \ref{partition} b). 
By construction, the cover $\mathcal{P}=\{I_\alpha\}_{\alpha<\omega_1}$, with $I_\alpha=[x_\alpha,x_{\alpha+1}]$, satisfies 
the first two claims.
Choose now for each $i=1,\dots,s$ a compact interval $J_i$ containing $c_i$ in its interior, and
set $J=J_1\times\cdots\times J_s$.
By Lemma \ref{pearl necklace lemma} applied to the cover $\mathcal{P}$, we get closed sets $B_\alpha\subset\RR^{n-s}$,
such that $B_\alpha$ 
contains $(I_\alpha)^{n-s}$ for all $\alpha>0$.
Therefore, $f(I_\alpha)\subset D_\alpha'\definition J\times B_\alpha$ for $\alpha>0$, and
the intersection 
$\partial D_\alpha'\cap f(I_\alpha)$ is equal to $\{f(x_\alpha),f(x_{\alpha+1})\}$.
By construction, $D_\alpha'\cap D_{\alpha+1}'=J\times\{(x_{\alpha+1},\dots,x_{\alpha+1})\}$.
Since $B_\alpha$ is a $(n-s)$-compact ball and
$f(I_\alpha)$ is contained in the hyperplane
$x_1=c_1,\dots,x_s=c_s$, one can find a subset $D_\alpha$ of $D_\alpha'$ which is diffeomorphic
to the compact ball in $\R^n$, that has the desired boundary properties (see Figure 3).
\hfill $\Box$
\begin{figure}[h]
  \begin{center}
  \epsfig{figure=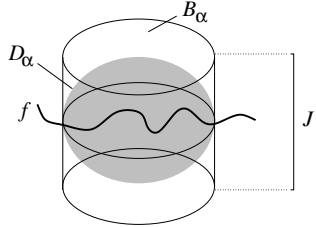, height=3cm}
  \caption{The subset $D_\alpha$ of $J\times B_\alpha$.}
\end{center}
\end{figure}
\begin{prop}
  \label{Schonfliess}
  Let $B$ be the closed unit ball in $\R^n$, and let
  $f,g:[0,1]\to B$ be continuous embeddings with $f^{-1}({\partial B})=g^{-1}({\partial B})=\{0,1\}$,
  $f(0)=g(0)$ and $f(1)=g(1)$. Then, there exists an isotopy $\phi_t:B\to B$ $(t\in[0,1])$ keeping $\partial B$ fixed,
  such that $\phi_0=id_B$ and $\phi_1\circ f=g$, if one of the following conditions holds:\\
  a) $n=2$,
  \\
  {b) $n\ge 4$ and $f,g$ are \CC{1} embeddings.}\\
  Suppose now that $f^{-1}({\partial B})=g^{-1}({\partial B})=\{0\}$ and $f(0)=g(0)$. Then, we have the same conclusion if:\\
  c) $n=2$, \\
  {d) $n\ge 3$ and $f,g$ are \CC{1} embeddings.}
\end{prop}
\noindent{\em Proof:}
Points a) and c) are direct consequences of Sch\"onfliess theorem. Point b) follows from Zeeman's unknotting theorem (see e.g.
\cite[Theorem 7.1]{RourkeSanderson}). Point d) follows from assertion b) for $n\ge 4$; for $n=3$, consider a good projection. 
\hfill $\Box$

%%%% sec 3

\section{Classification of long knots}\label{secres}

We are now ready to present our main results.
\begin{defi}\label{def direction}
The direction of a continuous function $f:\RR\to\RR^n$ is the vector
$D(f)=\left(\begin{array}{ccc}\delta_1(f)&\dots &\delta_n(f)\end{array}\right)^T$, where
$$
\delta_i(f)=\left\{\begin{array}{ll}
		1& \mbox{if $\pi_i\circ f$ is cofinal;}\\
		0& \mbox{if $\pi_i\circ f$ is bounded.}
		\end{array}\right.
$$
If $g:\RR\to\LL^n$ is continuous, let us define $\delta(g)=\left(\begin{array}{ccc}\delta_1(g)&\dots &\delta_n(g)\end{array}\right)^T$ as follows:
$$
\delta_i(g)=\left\{\begin{array}{lll}
		+1& \mbox{if $\pi_i\circ g$ is cofinal in $\LL_+$;}\\
		-1& \mbox{if $\pi_i\circ g$ is cofinal in $\LL_-$;}\\
		\phantom{-}0& \mbox{if $\pi_i\circ g$ is bounded.}
		\end{array}\right.
$$
The direction of a continuous fonction $f:\LL\to\LL^n$ is the $(n\times 2)$-matrix
$D(f)=\left(\begin{array}{cc}\delta(f|_{\LL_-})&\delta(f|_{\LL_+})\end{array}\right)$.
\\
A matrix $D$ is a $(Y,X)$-direction ($X=\RR,\LL$, $Y=\RR^n,\LL^n$) if $D=D(f)$ for some
continuous $f:X\to Y$.
\end{defi}
\begin{thm}\label{thm:invariant}
The direction is an invariant of $(\RR^n,\RR)$ and $(\LL^n,\LL)$-knots.
\end{thm}
\noindent{\em Proof:}
Let $f,g:\RR\to\RR^n$ be two continuous embeddings with $D(f)\neq D(g)$. Without loss of generality, it may be assumed that
$\pi_1\circ f$ is bounded and $\pi_1\circ g$ cofinal. Then, $f$ and $g$ are not equivalent. Indeed, consider an isotopy $\phi_t$ of $(\RR')^n$
such that $\phi_0=id$ and $\phi_1\circ f=g$; then, $\pi_1\circ\phi_t\circ f:\RR\to\RR'$ provides a homotopy between $\pi_1\circ f:\RR\to\RR\subset\RR'$
and $\pi_1\circ g:\RR\to\RR\subset\RR'$. Since $\pi_1\circ f$ is bounded, it follows from Lemma \ref{partition} a) that $\pi_1\circ f$ is homotopic to
a constant. On the other hand, Lemma \ref{partition} b) implies that $\pi_1\circ g$ is homotopic to the canonical inclusion $\RR\subset\RR'$. Therefore,
we would have a homotopy from a constant map to the inclusion $\RR\subset\RR'$. Since this inclusion is a homotopy equivalence, and since $\RR$
is not contractible (see e.g. \cite{Gauld}), such an homotopy does not exist. The proof for $(\LL^n,\LL)$-knots is very similar.
\hfill $\Box$

\subsection*{$(\RR^2,\RR)$ and $(\LL^2,\LL)$-knots}
\begin{thm}\label{thm1}
  The direction is a complete invariant for $(\RR^2,\RR)$-knots. Therefore, there are exactly $3$ classes of non-equivalent $(\RR^2,\RR)$-knots.
\end{thm}
\noindent{\em Proof:}
There are $4$ possible directions for a continuous map $f:\RR\to\RR^2$. If $f$ is an embedding, the direction $D=\verti{0}{0}$ is impossible.
(Otherwise, by Lemma \ref{partition} a), there would be some $z\in\RR$ such that $f$ is constant on $[z,\omega_1[$.) Therefore, we are left with three
possible directions, realized by the following $(\RR^2,\RR)$-knots: $f_{\verti{1}{0}}(x)=(0,x)$, $f_{\verti{0}{1}}(x)=(x,0)$ and $f_{\verti{1}{1}}(x)=(x,x)$.
By Theorem \ref{thm:invariant}, we just need to show that an $(\RR^2,\RR)$-knot $f$ with direction $D$ is equivalent to $f_D$, for
$D=\verti{1}{1},\verti{0}{1}$, and $\verti{1}{1}$.

Let us first assume that $D=\verti{1}{1}$, and let $\{I_\alpha=[x_\alpha,x_{\alpha+1}]\}_{\alpha<\omega_1}$ and $\{D_\alpha\}_{0<\alpha<\omega_1}$ be the covers
given by Lemma \ref{boxing lemma}.
Each $D_\alpha$ satisfies the hypotheses of Proposition \ref{Schonfliess} a),
with $f=f|_{I_\alpha}$ and $g=f_{\verti{1}{1}}|_{I_\alpha}$. Thus, we can find isotopies $(\phi_\alpha)_t$ of $D_\alpha$ $(rel \; \partial D_\alpha)$
such that $(\phi_\alpha)_0=id$ and $(\phi_\alpha)_1\circ f|_{I_\alpha}=f_{\verti{1}{1}}|_{I_\alpha}$. Then, consider a set $D_0\subset(\RR')^2$,
diffeomorphic to the compact ball, that contains $f([0,x_1[)\cup f_{\verti{1}{1}}([0,x_1[)$ in its interior, and
such that $D_0\cap D_1=\{(x_1,x_1)\}$. Then, $D_0$ satisfies the assumptions of Proposition \ref{Schonfliess} c), and there
in an ambient isotopy $(\phi_0)_t$ of $D_0$ $(rel \; \partial D_0)$ between $f|_{I_0}$ and $f_{\verti{1}{1}}|_{I_0}$. Extending the ambient
isotopies $(\phi_\alpha)_t$ $(\alpha<\omega_1)$ by the identity outside $\cup_{\alpha<\omega_1}D_\alpha$, we have proved that $f\sim f_{\verti{1}{1}}$.

Now, consider the case where $D=\verti{0}{1}$ (the case $D=\verti{1}{0}$ is similar). As before, take the cover
$\{D_\alpha\}_{0<\alpha<\omega_1}$ of $f([x_1,\omega_1[)$ given by Lemma \ref{boxing lemma}, and choose $D_0\subset(\RR')^2$,
diffeomorphic to the compact ball, that contains $f([0,x_1[)\cup (\{c_1\}\times[0,x_1[)$ in its interior, and
such that $D_0\cap D_1=\{(c_1,x_1)\}$. By Proposition \ref{Schonfliess} c), there
is an ambient isotopy of $D_0$ $(rel \; \partial D_0)$ that sends $f([0,x_1])$ on $\{c_1\}\times[0,x_1]$. Extending it by the identity
outside $D_0$, we have an ambient isotopy between $\text{\rm Im}f$ and $\{c_1\}\times\RR$.
It is now straightforward to define an isotopy between $\{c_1\}\times\RR$ and $\{0\}\times\RR=\text{\rm Im}f_{\verti{0}{1}}$. Therefore,
we can assume that $\text{\rm Im}f=\text{\rm Im}f_{\verti{0}{1}}=\{ (x,0)\, |\, x\in\RR\}$.
In that case, $\pi_1\circ f$ is an homeomorphism of $\RR$. By \cite[Corollary 2]{Gauld}, there is an isotopy $\gamma_t$
of $\RR$ (keeping $0$ fixed) such that $\gamma_0=id_\RR$ and $\gamma_1=\pi_1\circ f$.
Then, the isotopy $\phi_t$ of $(\RR')^2$ given by $\phi_t(x,y)=(\gamma_t(x),y)$ if $x\in\RR$ and $\phi_t(x,y)=(x,y)$ otherwise
satisfies $\phi_0=id$ and $\phi_1\circ f_D=f$.
\hfill $\Box$
\esp

The direction is almost a complete invariant for $(\LL^2,\LL)$-knots. It fails to be complete only
because some directions correspond to exactly
two non-equivalent knots. This phenomenon also appears in dimension $n\ge 2$, motivating the following definition.

\begin{defi}\label{double direction}
  A double direction is an $(\LL^n,\LL)$-direction with equal columns that contain exactly one $0$.
\end{defi}
\begin{thm}\label{corthm1}
   There are exactly $64$ non-equivalent $(\LL^2,\LL)$-knots.
\end{thm}
\noindent{\em Proof:}
Given $f:\LL\to \LL^2$ continuous, there are $9^2$ possible direction matrices $\left(\delta(f|_{\LL_-})\;\delta(f|_{\LL_+})\right)$. If $f$ is an
$(\LL^2,\LL)$-knot, the restrictions $f|_{\LL_-}$ and $f|_{\LL_+}$ are embeddings; by Lemma \ref{partition} a), 
this implies that the columns
$\delta(f|_{\LL_-})$ and $\delta(f|_{\LL_+})$ are non-zero. Furthermore, the following four directions are also forbidden:
$$
\scriptstyle
\left(\begin{array}{cc}+1&+1\\+1&+1\end{array}\right)\,,\;\;\left(\begin{array}{cc}+1&+1\\-1&-1\end{array}\right)\,,\;\;
\left(\begin{array}{cc}-1&-1\\+1&+1\end{array}\right)\,\;\mbox{and}\;\;\left(\begin{array}{cc}-1&-1\\-1&-1\end{array}\right)\,.
$$
(Indeed, let us consider a continuous map $f:\LL\to \LL^2$ with direction matrix 
% $\left(\begin{array}{cc}+1&+1\\+1&+1\end{array}\right)$.
$\verti{+1\, +1}{+1\, +1}$.
Since $\delta(f|_{\LL_+})=\verti{+1}{+1}$, it follows from Lemma \ref{partition} b) that the sets
$\{x\in\LL_+\;|\;\pi_1\circ f(x)=x\}$ and $\{x\in\LL_+\;|\;\pi_2\circ f(x)=x\}$ are
closed and cofinal. By Lemma \ref{closedcofinal}, so is their intersection $\{x\in\LL_+\;|\;f(x)=(x,x)\}$. Similarly,
$\delta(f|_{\LL_-})=\verti{+1}{+1}$ implies that $\{y\in\LL_-\;|\;f(y)=(-y,-y)\}$ is closed and cofinal.
By Lemma \ref{closedcofinal} again, the set $\{z\in\LL_+\;|\;f(z)=(z,z)=f(-z)\}$ is cofinal, so $f$ is not an embedding. The
other three cases are similar.) So, we are left with $(9-1)^2-4=60$ possible directions for an $(\LL^2,\LL)$-knot. It is easy to
exhibit a knot realizing each of these directions.

As in Theorem \ref{thm1}, one shows that
two knots with the same direction $D$ are equivalent, except if $D$ is a double direction.
To each double direction correspond exactly two classes of knots, as we shall see.
Since there are four double directions, namely
$$
\left(\begin{array}{rr}+1&+1\\0&0\end{array}\right)\,,\;\;\left(\begin{array}{rr}-1&-1\\0&0\end{array}\right)\,,\;\;
\left(\begin{array}{rr}0&0\\+1&+1\end{array}\right)\,\;\mbox{and}\;\;\left(\begin{array}{rr}0&0\\-1&-1\end{array}\right)\, ,
$$
the theorem follows.
Let $D$ be the first of these directions, and let $f:\LL\to\LL^2$
be an embedding with $f(x)=(x,1)$ for $x\ge 1$ and $f(x)=(-x,-1)$ for $x\le -1$ 
(as usual, we denote by $-x$ the point in $\LL_-$ corresponding to $x\in\LL_+$, and
vice versa).  Finally, let $g:\LL\to\LL^2$ be the knot given by $g(x)=f(-x)$. 
Clearly, $D=D(f)=D(g)$ and any $(\LL^2,\LL)$-knot with direction $D$ 
is equivalent to either $f$ or $g$. It remains to check that $f$ and $g$ are non-equivalent.
Indeed, let $\phi_t$ be an isotopy between $f$ and $g$, i.e.
$\phi_0=id$ and $\phi_1\circ f=g$. 
Since $\pi_2\circ\phi_t\circ f$ is a bounded continuous map $\LL\to\LL$, Lemma \ref{partition} a) 
implies that $\pi_2\circ\phi_t\circ f([x,\omega_1[)$ is a single point for $x$ large enough. 
Let us denote it by
$r_+(t)$, and similarly, let $r_-(t)$ be the element defined by $\pi_2\circ\phi_t\circ f(]-\omega_1,-x])$ for
$x$ large enough.
One checks that $r_+$ and $r_-$ are continuous. Since
$r_+(0)=r_-(1)=1$ and $r_+(1)=r_-(0)=-1$, there is some $t_0$ for which $r_+(t_0)=r_-(t_0)$. Then,
$\phi_{t_0}\circ f$ is not an embedding.
The other double directions are similarly treated.
\hfill $\Box$

\subsection*{Differentiable $(\RR^n,\RR)$ and $(\LL^n,\LL)$-knots for $n\ge 4$}

\begin{thm}\label{thm2}
  There are exactly $(2^n-1)$ non-equivalent differentiable $(\RR^n,\RR)$-knots if $n\ge 4$.
\end{thm}
\noindent{\em Proof:} Given $f$ an $(\RR^n,\RR)$-knot, at least one $\pi_i\circ f$ $(i=1,\dots,n)$ is cofinal;
hence, we have $2^n-1$ possible directions. Each of these directions $D$ can be realized by the $(\RR^n,\RR)$-knot $f_D$ given by $f_D(x)=x\cdot D^T$.
Now, we just need to prove that an $(\RR^n,\RR)$-knot $f$ with direction $D$ is equivalent to $f_D$.\\
By a permutation of the indices, it may be assumed that $\pi_i\circ f$ is bounded for $i=1,\dots,s$ and
cofinal for $i=s+1,\dots,n$. Consider the covers $\{[x_\alpha,x_{\alpha+1}]\}_{\alpha<\omega_1}$ of $\RR$ and
$\{D_\alpha\}_{0<\alpha<\omega_1}$ of $f([x_1,\omega_1[)$ given by Lemma \ref{boxing lemma}. For $i=1,\dots,s$, take an isotopy $\phi^i_t$ of $\RR'$
between $\pi_i\circ f(x_1)=c_i$ and $0$; the isotopy $\phi_t=(\phi^1_t,\dots,\phi^s_t,id,\dots,id)$ then sends $f([x_1,\omega_1[)$ on $\{0\}\times\RR^{n-s}$.
Using Proposition \ref{Schonfliess} d), we may assume that $f=f_D$ on $[0,x_1]$. By Proposition \ref{Schonfliess} b), we have an isotopy $(\phi_\alpha)_t$
in $D_\alpha$ between $f|_{[x_\alpha,x_{\alpha+1}]}$ and $f_D|_{[x_\alpha,x_{\alpha+1}]}$ for all $0<\alpha<\omega_1$. Extending this isotopy by the identity
outside $\cup_\alpha D_\alpha$, it follows that $f\sim f_D$.
\hfill $\Box$

\begin{thm}\label{corthm2}
 There are exactly $(3^n-1)^2+(n-2)2^{n-1}$
 non-equivalent differentiable $(\LL^n,\LL)$-knots if $n\ge 4$.
\end{thm}
\noindent{\em Proof:} We have $(3^n)^2$ possible directions for a continuous map $f:\LL\to\LL^n$. If $f$ is an embedding,
$\delta(f|_{\LL_-})$ and $\delta(f|_{\LL_+})$ are non-zero, giving $(3^n-1)^2$ directions. Furthermore, a direction matrix with
identical columns and no zero coefficient cannot be realized by an embedding (the argument is similar to the case $n=2$).
Therefore, we are left with $(3^n-1)^2-2^n$ directions, which can all clearly be realized by $(\LL^n,\LL)$-knots.
Finally, to any $(\LL^n,\LL)$-direction corresponds exactly one class of knots, except if
it is a double direction; in this case, there are exactly two classes of knots with this direction.
We thus have $(3^n-1)^2 - 2^n+n 2^{n-1}=(3^n-1)^2+(n-2)2^{n-1}$ different classes of knots.
\hfill $\Box$

\subsection*{Differentiable $(\RR^3,\RR)$ and $(\LL^3,\LL)$-knots}
\begin{thm} \label{thm3}
  There are $2^3-1=7$ non-equivalent differentiable $(\RR^3,\RR)$-knots.
\end{thm}
\noindent We shall need these two lemmas:
\begin{lemma}\label{monotone}
  Let $f$ be a differentiable $(\RR^3,\RR)$-knot such that $\pi_i\circ f$ is cofinal for $i=1,2,3$, and let
  $\mathcal{P}^3\definition\{[x_\alpha,x_{\alpha+1}]^3\}_{0<\alpha<\omega_1}$ with $x_\alpha$ as in Lemma \ref{boxing lemma}.
  Then, for all but finitely many $\alpha$, there is an index $i$ with $\pi_i\circ f$ monotone on $[x_\alpha,x_{\alpha+1}]$.
\end{lemma}
\begin{lemma} \label{monotone2}
  Let $f:[0,1]\to[0,1]^3$ be a differentiable embedding, and let $B$ be a compact $3$-ball such that $B\supset[0,1]^3$ and
  $\partial B\cap[0,1]^3=\{(0,0,0),(1,1,1)\}$. Suppose that $f(0)=(0,0,0)$, $f(1)=(1,1,1)$, and that $\pi_i\circ f$ is monotone for at least one $i$;
  then, there is an isotopy $\phi_t$ of $B$, which is the identity on $\partial B$, such that $\phi_0=id$ and $\phi_1\circ f(x)=(x,x,x)$.
\end{lemma}
\noindent{\em Proof of Lemma \ref{monotone}:}
Otherwise, let $\{\beta_m\}_{m<\omega}$ be an increasing sequence of ordinals such that $\pi_i\circ f$ is not monotone on
$I_{\beta_m}=[x_{\beta_m},x_{\beta_m+1}]$ for $i=1,2,3$, and let us denote by $\beta$ the limit of this sequence. For $m<\omega$ and
$i=1,2,3$, choose $x_{i,m},y_{i,m}\in I_{\beta_m}$ such that
$\pi_i\circ f$ is decreasing in a neighborhood of $x_{i,m}$ and increasing in a neighborhood of $y_{i,m}$.
In any chart containing $[0,x_\beta]$,
$(\pi_i\circ f)'(x_{i,m})\le 0$ and $(\pi_i\circ f)'(y_{i,m})\ge 0$ for any $m<\omega$ and $i=1,2,3$.
By construction, $\displaystyle\lim_{m\to\infty}x_{i,m}=\lim_{m\to\infty}y_{j,m}\definition u$ for any $i,j$.
By continuity,
$(\pi_i\circ f)'(u)=0$ for $i=1,2,3$, contradicting the $1$-regularity of $f$.
\hfill $\Box$
\esp
\noindent{\em Proof of Lemma \ref{monotone2}:}
Apply a descending curve argument.
\hfill $\Box$
\esp
\noindent{\em Proof of Theorem \ref{thm3}:}
Clearly, there are $7$ possible $(\RR^3,\RR)$-directions. Let us prove that two knots $f,g$ with $D(f)=D(g)$ are equivalent.\\
First, consider the case where $\pi_i\circ f$ is bounded for some index $i$ (let us say
that $i=1$). By Lemma \ref{partition} a), there is an $x_1$ in $\RR$ such that $\pi_1\circ f(x)=c$ for all $x\ge x_1$. Let $D_0$ be a compact
$3$-ball in $(\RR')^3$ such that $f([0,x_1])\subset D_0$, $f^{-1}(\partial D_0)=\{x_1\}$ and $f([0,x_1])\cap D_0=\{f(x_1)\}$. By Proposition
\ref{Schonfliess} d), $f$ is isotopic to an $(\RR^2,\RR)$-knot in $\{c\}\times\RR^2$. We can conclude with Theorem \ref{thm1}.\\
Now, let us assume that $\pi_i\circ f$ is cofinal for $i=1,2,3$. By
Lemmas \ref{monotone} and \ref{monotone2}, there is an $x$ in $\RR$ such that $f|_{[x,\omega_1[}$ is isotopic (in $[x,\omega_1[^3$) to the
knot $f_D$ given by $f_D(x)=(x,x,x)$. We then proceed as before for $f|_{[0,x]}$.
\hfill $\Box$

Let us now turn to $(\LL^3,\LL)$-knots. We will assume that $\mathbb{S}^3$
and $\LL^3$ have a fixed orientation.
Let $k$ be an oriented differentiable $(\mathbb{S}^3,\mathbb{S}^1)$-knot, and let $D$ be an $(\LL^3,\LL)$-direction.
If $D$ is not a double direction,
there is clearly a
unique equivalence class of `unknotted' differentiable $(\LL^3,\LL)$-knot with direction $D$; as before, let us denote a representant by $f_D$.
If $D$ is a double direction, there are exactly two such `unknotted' $(\LL^3,\LL)$-knots with direction $D$; 
to simplify the notation we shall 
abusively denote both by $f_D$.

\begin{figure}[h]
  \begin{center}
    \epsfig{figure=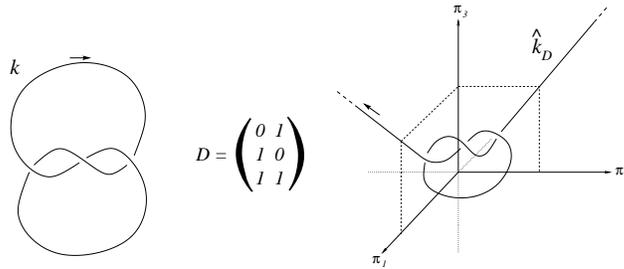, height=3.5cm}
    \caption{The trefoil knot and its $D$-associate.}
    \label{fig:associate}
  \end{center}
\end{figure}

\begin{defi}
  A differentiable $(\LL^3,\LL)$-knot $\wh{k}_D$ equivalent to the oriented connected sum $f_D\;\#\;k$ is called a $D$-associate of $k$.
  (See Figure \ref{fig:associate} for an example.)
\end{defi}
It is easy to show that, up to the subtleties due to the double directions, 
the equivalence class of $\wh{k}_D$ only depends on the equivalence class of the oriented
$(\mathbb{S}^3,\mathbb{S}^1)$-knot $k$ and on the direction $D$: just follow the classical proof that the connected sum of
oriented differentiable $(\mathbb{S}^3,\mathbb{S}^1)$-knots is well defined.
Counting the directions as in Theorem \ref{corthm2}, we find immediately:
\begin{prop}\label{prop:ass}
  A differentiable $(\mathbb{S}^3,\mathbb{S}^1)$-knot has exactly $680$
  non-equivalent associates.\hfill $\Box$
\end{prop}
Moreover, the classification of smooth $(\LL^3,\LL)$-knots with direction $D$ is equivalent to the classification of
differentiable oriented $(\mathbb{S}^3,\mathbb{S}^1)$-knot. In other words:
\begin{thm}\label{thm4}
  A differentiable $(\LL^3,\LL)$-knot is the associate of a unique type of differentiable oriented $(\mathbb{S}^3,\mathbb{S}^1)$-knot.
\end{thm}
\noindent{\em Proof:} Let $f$ be a differentiable $(\LL^3,\LL)$-knot with direction $D$. Using the same argument as in Theorem \ref{thm3},
we see that for some $x_1\in\LL_+$, there is an isotopy $\varphi_t$ of $X=\LL^3\backslash]-x_1,x_1[^3$, keeping $\partial X$ fixed, such
that $\varphi_0=id_X$ and $\varphi_1\circ f|_{f^{-1}(X)}=f_D|_{f_D^{-1}(X)}$. Let $k$ be the smooth oriented $(\mathbb{S}^3,\mathbb{S}^1)$-knot
obtained by attaching both ends of $f|_{f^{-1}([-x_1,x_1]^3)}$  with an unknotted arc in $X$. By construction, $f$ is equivalent to $f_D\;\#\;k$, so
$f$ is the $D$-associate of $k$.
\\
Now, let $k,k'$ be two differentiable $(\mathbb{S}^3,\mathbb{S}^1)$-knots such that there is an ambient isotopy $\phi_t$ of $\LL^3$ with
$\phi_1\circ\wh{k}_D=\wh{k}'_D$. If $U_\alpha$ denotes the open cube $]-\alpha,\alpha[^3$ in $\LL^3$, then
$\{U_\alpha\,|\,\alpha<\omega_1\}$ is a {\em canonical sequence} in the sense of \cite[p. 147]{Gauld}.
By the proposition on the same page, the set
$$
  F_\phi\definition
  \{\alpha<\omega_1\,|\, \phi_t({\overline U_\alpha}\backslash U_\alpha)=
  {\overline U_\alpha}\backslash U_\alpha\, \forall t\}
$$
is cofinal in $W(\omega_1)$.
(The proposition is stated for $\omega$-bounded $2$-surfaces, but its proof shows
that it also holds for $\RR^n$ and $\LL^n$.) Take $\alpha$ large enough such that $\wh{k}_D|_{\LL^3\backslash U_\alpha}$ and
$\wh{k}'_D|_{\LL^3\backslash U_\alpha}$ are equal to $f_D|_{\LL^3\backslash U_\alpha}$; let us denote by $\{P_1,P_2\}$ the two points of
$\text{\rm Im}\wh{k}_D\cap({\overline U_\alpha}\backslash U_\alpha)=\text{\rm Im}\wh{k}'_D\cap({\overline U_\alpha}\backslash U_\alpha)$.
For all $t$, $\phi_t|:{\overline U_\alpha}\backslash U_\alpha\to{\overline U_\alpha}\backslash U_\alpha$ is a homeomorphism
isotopic $(rel\;\{P_1,P_2\})$ to the identity. Therefore, $\phi_t|$ can be extended to a homeomorphism $\widetilde{\phi}_t:B\to B$,
where $B$ is a compact ball containing ${\overline U_\alpha}$, with
$\widetilde{\phi}_t$ keeping $\partial B$ fixed. Extending $\widetilde{\phi}_t$ with the identity on $\mathbb{S}^3\backslash B$,
we get an isotopy between $k$ and $k'$.
\hfill $\Box$

%%%%%%%%

\footnotesize

\bibliographystyle{plain}
\bibliography{/home/baillif/macros/biblio}

\vskip1cm
\noindent
Mathieu Baillif \\
Section de Math\'ematiques\\
2-4 rue du Li\`evre\\
1211 Gen\`eve 24\\
Switzerland\\
{\tt Mathieu.Baillif@math.unige.ch}

\vskip1cm
\noindent
David Cimasoni\\
Section de Math\'ematiques\\
2-4 rue du Li\`evre\\
1211 Gen\`eve 24\\
Switzerland\\
{\tt David.Cimasoni@math.unige.ch}

\end{document}